# A fuzzy framework for unraveling experiential tourism value

Eddy SORIA[a,1], Aïda VALLS[b], Jorge Luis MARIÑO[c] and Ana Beatriz HERNÁNDEZ[a]
[a] *Department of Business Management*
*Universitat Rovira i Virgili, Catalonia, Spain*
[b] *Department of Computer Science and Mathematics*
*Universitat Rovira i Virgili, Catalonia, Spain*
[c] *Department of Tourism*
*Universidad de Oriente, Cuba*
ORCiD ID: Eddy Soria https://orcid.org/0000-0001-9999-5850, Aïda Valls https://orcid.org/0000-0003-3616-7809, Jorge Luis Mariño https://orcid.org/0000-0003-4663-9689 and Ana Beatriz Hernández https://orcid.org/0000-0002-8110-9328

**Abstract.** There is currently a growing interest in the study of sustainable tourism experiences. Although the main dimensions of tourism experiences are known, the assessment of attractions as a key component in sustainable tourism design has seen little research. This gap stems from two primary challenges: firstly, the complexity of each tourism destination, which makes use of heterogeneous scales of evaluation; and secondly, the subjective nature of tourism experiences that requires qualitative measurement scales and assessment. To overcome these shortcomings, a new methodology based on fuzzy numbers is proposed, consisting of a scale normalization stage, followed by an aggregation of the fuzzy numbers that can be integrated into a Geographical Information System. This technique is tested on the quantification of the tourism value of attractions for sustainable experience in Santiago de Cuba.

## 1. Introduction

The tourism sector plays a significant role in promoting employment and boosting the economy in various regions, thereby contributing to poverty reduction and the socio-economic improvement of several countries. This sector is comprised of numerous components such as tourist attractions, facilities, infrastructure, accessibility, tourism demand and competition, among others, which interact to shape the so-called “tourism experience” [1]. Particularly, tourist attractions are considered the primary requirement of tourism [2], and constitute the fundamental element for the creation of experiences in tourist destinations.

The term “experience” denotes a specific type of tourism product. A tourism product refers to the tangible and intangible components offered to tourists, such as attractions that motivate people to travel, and often includes accommodation, transportation, and excursions. The tourism experience, on the other hand, is the personal and emotional outcome that tourists derive from their relationship with these products, encompassing their interactions, satisfaction and overall impression of the trip.

[1] Corresponding Author: Eddy Soria Leyva, eddy.soria@urv.cat

Both researchers and destination stakeholders are constantly asking themselves how to sell tourism experiences that increase the number of visitors without negatively impacting the environment [3]. As a result, some attempts have been made to align the design of experiences with the principles of sustainable development, producing what is referred to as a sustainable tourism experience [4]. Thus, a sustainable tourism experience can be defined as a tourism product capable of providing the tourist with profound emotions and memories, arising from stimuli that act on their overall perception and motivate them to contribute to the sustainability of the destination.

The process of designing sustainable tourism experiences is not easy. A significant part of the complexity lies in measuring and quantifying the dynamically intervening factors at each destination and as a result, universal techniques are not usually developed. On the contrary, the measurement indicators -usually known as Tourism Value Index (TVI)- are adapted to the characteristics of each destination, without establishing a standardization or normalization process, making hardly comparable the studies carried out in different tourism destinations, or even in the same tourist destination at different times.

In addition, tourism experiences are influenced by many subjective factors that define tourists' motivations and desires. Therefore, measuring these factors relevant for the design of tourism experiences requires the modeling of imprecise subjective evaluations. However, a thorough review of the literature has revealed that existing studies have not focused directly on addressing the subjectivity inherent in the calculation of TVI, nor have they attempted to adapt TVI to the new context of sustainable tourism experiences. Instead, it has traditionally been addressed in the classic framework of tourism attractions. In fact, the very concept of sustainable tourism experiences has only been loosely used without a proper definition or empirical operationalization [5].

The goal of this research is to address these shortcomings by introducing a new technique based on fuzzy numbers that can quantify the different indicators of the tourism value of attractions, transform them into a common scale and aggregate them. Therefore, the main contribution of this paper is the proposal of the Fuzzy Experiential Tourism Value (FTV) index, which provides a more comprehensive and flexible measure for supporting the design of sustainable tourism experiences than traditional methods of calculating TVI. While conventional TVI calculations rely on crisp numerical data, the proposed FTV integrates Triangular Fuzzy Numbers along with a normalization process. This approach generates adaptable tourism value estimates better suited for analytical decision-making in sustainable tourism experiences.

Section 2 gives some hints about the studies on sustainable tourism experiences and presents the most common definitions of Tourism Value Indexes in the literature. Section 3 explains the proposed methodology for constructing the new fuzzy approach to the tourism value index. First, we introduce a data normalization technique tailored for uncertain contexts, called Fuzzy Linear Rescaling, which is an extension of Linear Rescaling [1] using Triangular Fuzzy Numbers (TFNs). Secondly, the Fuzzy Experiential Tourism Value (FTV) is formulated based on the Fuzzy Linear Rescaling. In section 4 the proposed FTV index is evaluated in a real-world case study to obtain the FTV of the tourist attractions at Santiago de Cuba, and the results are integrated into a Geographic Information System (GIS) for spatial analysis and visualization. Section 5 gives the conclusions and outlines some future research open directions.

## 2. Related Work

**Sustainable tourism experiences.** The study of tourism experiences has historically played a significant role, as evidenced by its presence in at least three of the five most cited publications within tourism, leisure, and hospitality journals [6, 7]. Concurrently, since 2008, there has been a marked rise in the focus on sustainable tourism, highlighting the field's growing importance as indicated by [6]. However, according to [5], at the beginning of the current decade, there had not yet been a clear conceptualization and operational definition of the concept of sustainable tourism experiences.

Sustainable development is the harmonization of economic growth, social inclusion, and environmental protection to ensure long-term prosperity and well-being for current and future generations. By incorporating sustainability into the design of experiences, the development of new sustainable tourism experiences can be shaped. For example, some eco-lodges in Costa Rica are designed to minimize environmental impact [8] using renewable energy sources, water conservation measures, and sustainable building materials. Through these practices, tourists experience a profound connection with nature, gain insight into the importance of environmental conservation, and are motivated to contribute to the sustainability of the ecosystem, thereby embodying the essence of sustainable tourism experiences.

**Tourism Value Index**. Tourism attractions are the foundational elements of tourism products [2]. Consequently, they also play a pivotal role in the design process of sustainable tourism experiences, as there can be no experiences for tourists without attractions to visit. In general, the assessment of attractions is usually conducted through the calculation of the traditional Tourism Value Index (TVI) which provides a crisp measure of the level of attractiveness capable of generating demand from the attractions. TVI is a synthetic indicator that merges the factors involved in the tourist experience of a given destination. This indicator has been previously defined in different ways in the literature to rank tourism attractions, as we can see in **Table 1**. However, readapting TVI specifically to sustainable tourism experiences with subjective criteria has not been previously considered.

**Table 1.** Traditional approaches for calculating TVI. Notation: $(S_i)$ Score Value; $(\omega_i)$ Weight of the i-th criterion; (n) Number of criteria; (m) Number of experts; and for hierarchical criteria: $(S_{ji})$ Score Value for the i-th sub-criterion in the j-th criteria set; $(\omega_{ji})$ Weight of the i,j sub-criterion.

| Reference | Definition | |
|---|---|---|
| [9] | $TVI = \sum_{i=1}^{n} S_i * \omega_i$ | (1) |
| [10] | $TVI = \frac{\sum_{i=1}^{n} S_i * (1+\omega_i)}{\sum_{i=1}^{n} (1+\omega_i)}$ | (2) |
| [11] | $TVI = \sum_{i=1}^{n} S_i * \omega_i$ | (3) |
| [12] | $TVI = \sum_{i=1}^{n} \omega_j (\omega_{ji} * S_{ji})$ | (4) |
| [13] | $TVI = \frac{5}{\mathrm{n}} \sum_{\mathrm{j}=1}^{m} \sum_{\mathrm{i}=1}^{n} \left( \omega_{\mathrm{i}} \frac{(\mathrm{S_{ij}} - \min(\mathrm{S_{ij}}))}{\max(\mathrm{S_{ij}}) - \min(\mathrm{S_{ij}})} \right)$ | (5) |

**Data Normalization**. In the TVI definitions the values of the input criteria must be all in the same measurement scale, although they are usually collected from indicators that use different scales. It is common to use different scales of measurement for the indicators, depending on their meaning. For example, we can use a scale from 0 to 100 to indicate the percentage of accomplishment of a factor, a Likert scale from 1 to 5 to measure the tourist satisfaction, and a scale from -5 to 0 for the accessibility difficulties to an attraction. From a statistical standpoint, data normalization refers to the process of adjusting values measured on different scales to a common scale. There are several techniques for normalizing data, including standardization (z-score), rescaling, decimal scaling normalization, and normalization with respect to the maximum. In this research, we will use normalization based on rescaling to the interval [m, M], since it maintains the relationships between the original data values [1]:

**Definition 1.** Let $\mathrm{A} = \{\mathrm{a} \in \mathbb{R} \mid \mathrm{x} \leq \mathrm{a} \leq \mathrm{y};\ \mathrm{x} \neq \mathrm{y}\}$ be a nonempty set. Linear Rescaling is a mapping LRE: $[\mathrm{x}, \mathrm{y}] \rightarrow [\mathrm{m}, \mathrm{M}]$ defined by:

$$LRE\ (a) = M - (M - m)\left(\frac{y - a}{y - x}\right) \quad (6)$$

## 3. Proposed Method

**Definition 2.** Let $\tilde{a} = \left(a^{(1)}, a^{(2)}, a^{(3)}\right)$ be a Triangular Fuzzy Number (TFN) where $a^{(i)} \in [x, y]$ with $x, \mathrm{y} \in \mathbb{R}$ and $\mathrm{x} \neq \mathrm{y}$, and with $a^{(1)} \leq a^{(2)} \leq a^{(3)}$. Then, a Triangular Fuzzy Rescaling (TFR) is a mapping $TFR: \tilde{a} \rightarrow [m, M]$ defined by:

$$TFR\ (\tilde{a}) = \left(M - (M - m)\left(\frac{y - a^{(1)}}{y - x}\right), M - (M - m)\left(\frac{y - a^{(2)}}{y - x}\right), M - (M - m)\left(\frac{y - a^{(3)}}{y - x}\right)\right) \quad (7)$$

, where $m$ and $M$ are two parameters such that $m \in \mathbb{R}, M \in \mathbb{R}, m < M$.

Triangular Fuzzy Rescaling (TFR) is a generalized form of Linear Rescaling, extended to the realm of TFNs.

The Fuzzy Experiential Tourism Value is computed by taking a weighted average of the normalized Triangular Fuzzy Numbers as follows.

**Definition 3:** Fuzzy Experiential Tourism Value

Let $\tilde{a} = (\tilde{a}_1, \tilde{a}_2, \dots \tilde{a}_n)$ be a set of $n$ fuzzy numbers $\Psi$ where $\tilde{a}_i$ is a Triangular Fuzzy Number that satisfies $\tilde{a}_i \in [x_i, y_i], [x_i, y_i] \in IN$ ; and let $IN = \{[X, Y] \mid X \in \mathbb{R}, Y \in \mathbb{R}, X \neq Y\}$ be a non-empty set of interval measurement scales. The Fuzzy Experiential Tourism Value ($FTV$) is a function $FTV: \tilde{a} \rightarrow [m, M]$ that aggregates factors as defined in Eq. 8. The function is associated with an $n$-dimensional weighting vector $\omega$, where $\omega_i \in [0,1]$ and $\sum_{i=1}^{n} \omega_i = 1$.

$$FTV\ (\tilde{a}) = \begin{pmatrix} \sum_{i=1}^{n} \omega_i \left(M - (M - m)\left(\frac{y_i - a_i{}^{(1)}}{y_i - x_i}\right)\right), \\ \sum_{i=1}^{n} \omega_i \left(M - (M - m)\left(\frac{y_i - a_i{}^{(2)}}{y_i - x_i}\right)\right), \\ \sum_{i=1}^{n} \omega_i \left(M - (M - m)\left(\frac{y_i - a_i{}^{(3)}}{y_i - x_i}\right)\right) \end{pmatrix} \quad (8)$$

## 4. Case Study

### *4.1. Data*

The historical center of Santiago de Cuba was selected as a case study due to its high relevance for tourism economy in this city [1]. Interviews were conducted with residents, tourism sector workers, and tourists to understand their perspectives and determine the set of attractions to study, which consists on a set of 245 tourist attractions.

A total of 20 factors were included in this study, which are listed in Table 2. They were selected after revising the literature and considering the characteristics of Santiago de Cuba. Specific measurements were carried out on-site by an expert coder and reviewed by another [14]. The weightings were assigned based on the Allocation Hypothesis Method according to their importance for development, following [15], and normalized relative to the total so that the weighting vector summed to 1.

**Table 2.** Factors for measuring FTV.

| N | Factors | References | N | Factors | References |
|---|---|---|---|---|---|
| **1** | Attraction Condition | [10, 16, 17] | **11** | Customization Level | [18] |
| **2** | Attraction Info | [10, 16, 19] | **12** | Heritage Value | [5, 16] |
| **3** | Visitor Type | [10, 17] | **13** | Educational Value | [16, 18, 20] |
| **4** | Visit Date | [10, 17] | **14** | Sensory Impact | [18] |
| **5** | Historical Origin | [10, 16] | **15** | Economic Impact | [16, 21] |
| **6** | Residents' View | [5, 10] | **16** | Social Impact | [21, 22] |
| **7** | Tourists' View | [5, 10, 16, 17] | **17** | Authenticity | [18, 23] |
| **8** | Tourist Engagement | [5, 16, 18] | **18** | Environmental Impact | [5, 10, 16] |
| **9** | Stakeholders' Engagement | [5] | **19** | Cultural Commercialization | [10, 22] |
| **10** | Enjoyment Level | [16, 17, 20] | **20** | Social Vulnerability | [10] |

After calculating the FTV of each attraction, the obtained fuzzy number was defuzzified following with the method proposed by [24], to express it on a scale [0-100] in order to distinguish three levels of attractiveness (Low: 0-33; Medium: 34-66; and High: 67-100). The attractions in High category (FTV>=67) were selected for their inspection by the destination managers. Additionally, the FTV index was used to aggregate the values of all attractions per each factor in order to show the performance of the different criteria. The geographical QGIS software was then used to codify the fuzzy experiential tourism value measurements for each attraction in a geographic information system. This software was used to create heat maps based on the density of the experiential tourism value.

### *4.2. Results*

In this study area, 245 tourist attractions were studied. First, the factors listed in Table 2 were analyzed by applying the FTV index to aggregate the data of the 245 together per criterion (see Table 3 and Figure 1). As result, the factors with the highest levels of FTV normalized to the interval [0, 100] were: “Cultural Commercialization” (78,98,100), “Social Vulnerability” (78, 98, 100), “Environmental Impact” (76, 97, 100), "Visitor Type" (78, 98, 99), "Attraction Info" (72, 92, 97), "Attraction Condition" (71, 91, 97) and "Visit Date" (68, 88, 91).

**Table 3.** Assessment of touristic and sustainability factors using FTV

| Factors | x | y | $\omega_i$ | Non-normalized TFN | Normalized TFN |
|---|---|---|---|---|---|
| Attraction Condition | 0.00 | 5.00 | 0.054 | (3.54, 4.54, 4.87) | (71,91,97) |
| Attraction Info | 0.00 | 5.00 | 0.032 | (3.60, 4.60, 4.85) | (72,92,97) |
| Visitor Type | 0.00 | 5.00 | 0.043 | (3.90, 4.90, 4.97) | (78,98,99) |
| Visit Date | 0.00 | 5.00 | 0.038 | (3.42, 4.42, 4.56) | (68,88,91) |
| Historical Origin | 1.00 | 5.00 | 0.054 | (1.19, 1.87, 2.86) | (5,21,46) |
| Residents' View | 0.00 | 5.00 | 0.052 | (3.11, 4.11, 4.73) | (62,82,95) |
| Tourists' View | 0.00 | 5.00 | 0.055 | (3.06, 4.06, 4.73) | (61,82,95) |
| Tourist Engagement | 1.00 | 5.00 | 0.051 | (2.92, 3.92, 4.61) | (48,74,91) |
| Stakeholders' Engagement | 1.00 | 5.00 | 0.050 | (3.13, 4.13, 4.90) | (53,78,98) |
| Enjoyment Level | 1.00 | 5.00 | 0.052 | (2.44, 3.44, 4.12) | (36,61,79) |
| Customization Level | 1.00 | 5.00 | 0.049 | (1.97, 2.97, 3.92) | (24,50,74) |
| Heritage Value | 1.00 | 5.00 | 0.040 | (2.96, 3.96, 4.69) | (49,74,93) |
| Educational Value | 1.00 | 5.00 | 0.042 | (2.40, 3.40, 4.37) | (35,60,84) |
| Sensory Impact | 1.00 | 5.00 | 0.055 | (2.66, 3.66, 4.49) | (42,67,88) |
| Economic Impact | 1.00 | 5.00 | 0.053 | (2.35, 3.35, 4.31) | (34,59,82) |
| Social Impact | 1.00 | 5.00 | 0.053 | (3.28, 4.28, 4.87) | (57,82,97) |
| Authenticity | 1.00 | 5.00 | 0.057 | (3.01, 4.01, 4.60) | (50,75,90) |
| Environmental Impact | -5.00 | 0.00 | 0.057 | (-1.18, -0.18, -0.02) | (76,97,100) |
| Cultural Commercialization | -5.00 | 0.00 | 0.055 | (-1.10, -0.10, 0.00) | (78,98,100) |
| Social Vulnerability | -5.00 | 0.00 | 0.056 | (-1.10, -0.10, 0.00) | (78,98,100) |

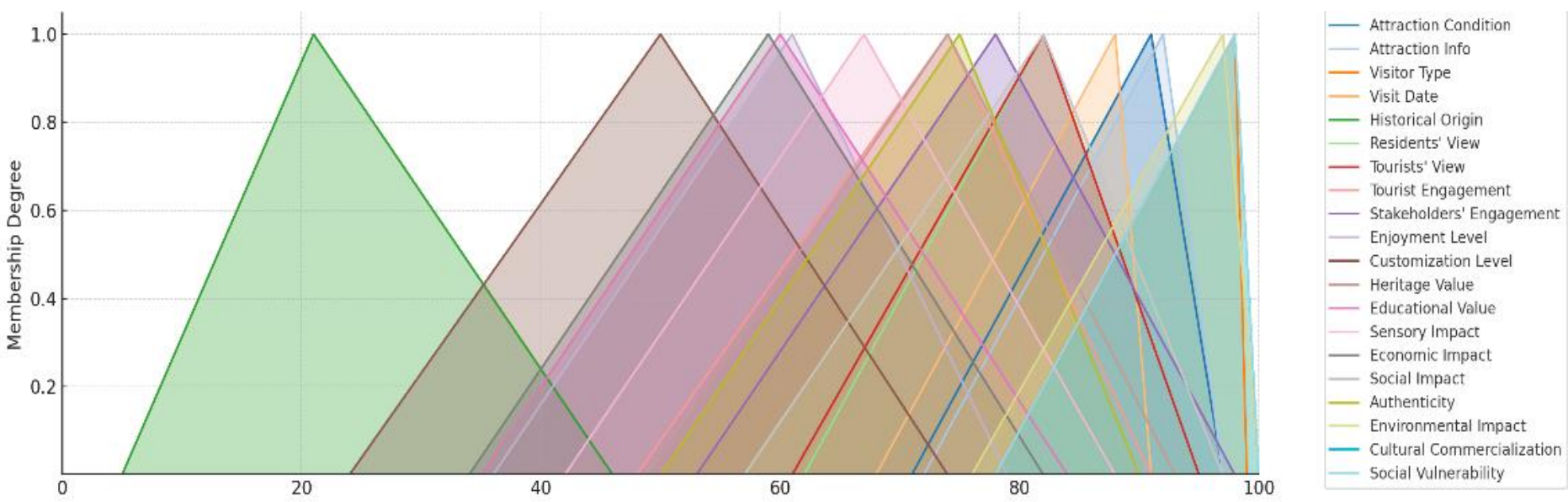

**Figure 1**. Normalized TFN Membership Functions

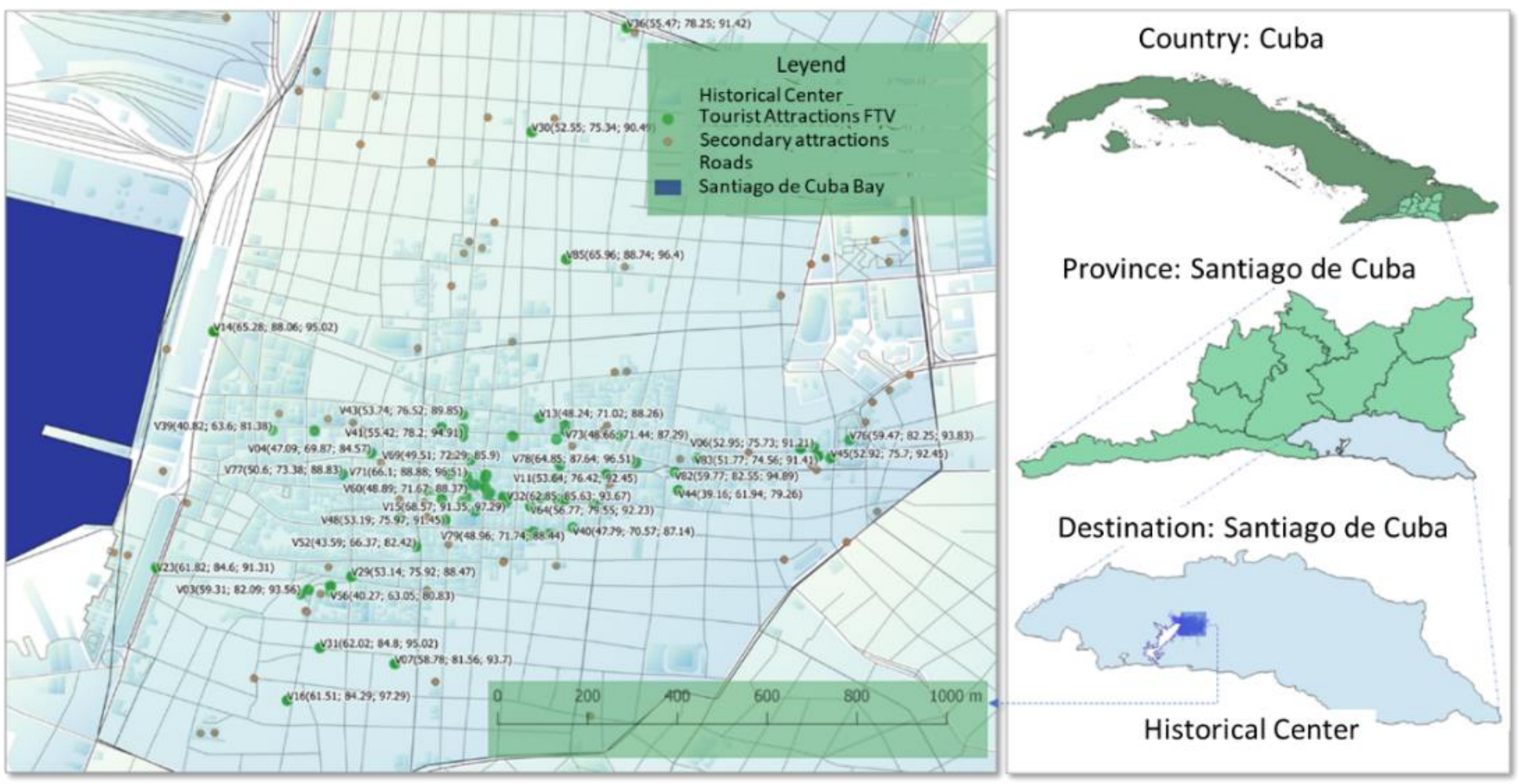

**Figure 2**. FTV distribution

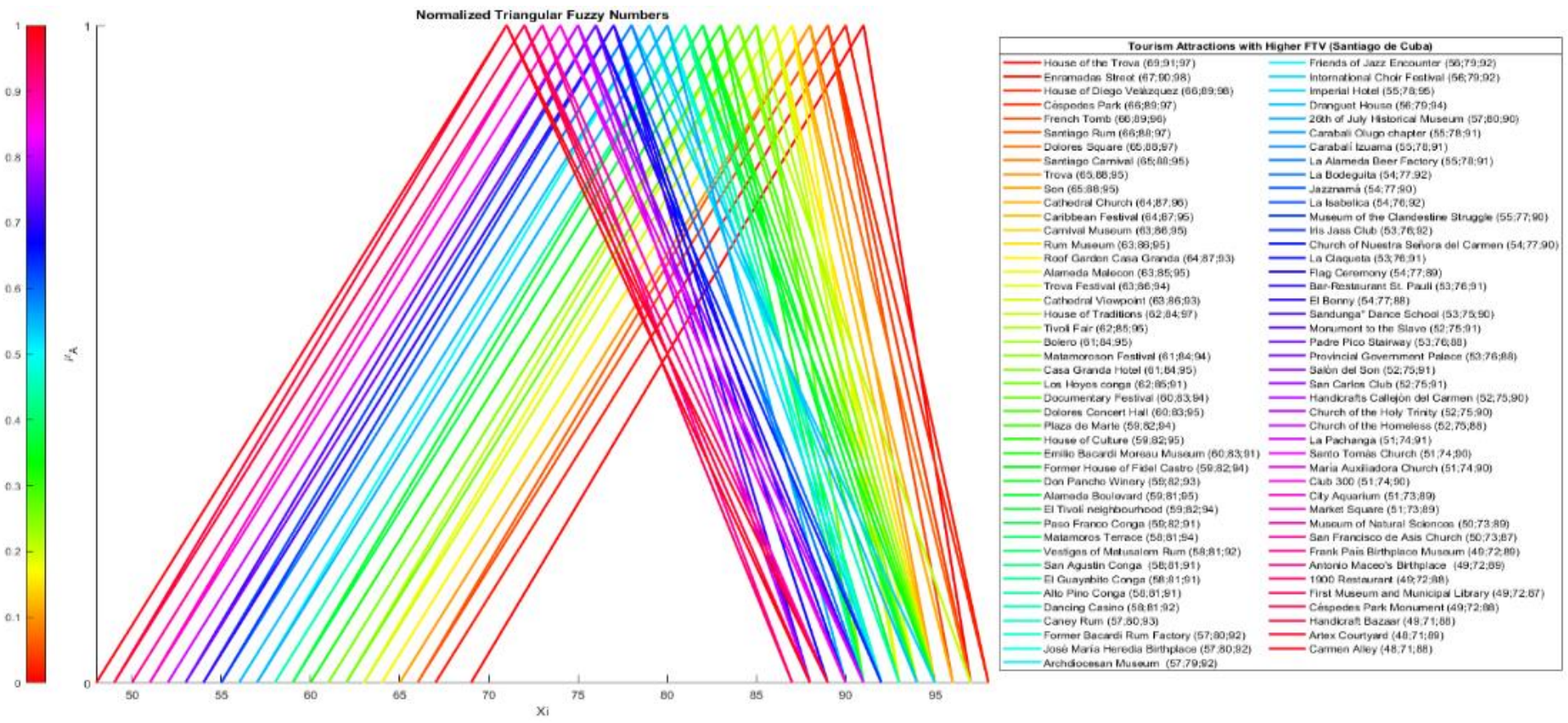

**Figure 3**. Filtered attractions with higher Fuzzy Experiential Tourism Value (FTV)

Secondly, the FTV index was applied to the attractions. The results were mapped in QGIS software and are visualized in Figure 2. From them, 87 attractions with the most representative FTVs were selected. Their final FTV fuzzy value is given in Figure 3.

The attractions with the highest FTV scores were: “House of the Trova” (FTV: 68.57, 91.35, 97.29), “Enramadas Street” (FTV: 67.32, 90.1, 97.73), “House of Diego Velázquez” (FTV: 65.96, 88.74, 97.56), “Céspedes Park” (FTV: 66.1, 88.88, 96.51), and “The French Tomb Society La Charité d’Orient” (FTV: 65.96, 88.74, 96.4). Next, the concentration degree of the 87 representative attractions distributed throughout the historic center was analyzed in QGIS to identify the areas with the highest tourism value and potential for sustainable tourism experiences. Six hot spots were visually identified through spatial analysis of the density of attractions (Figure 4): the area with the highest concentration was "Parque Céspedes", followed by: "Loma del Intendente", "Calle Heredia", "Calle Enramadas", "Plaza Dolores" and "Plaza de Marte". Sustainable tourism experiences can be developed in these areas.

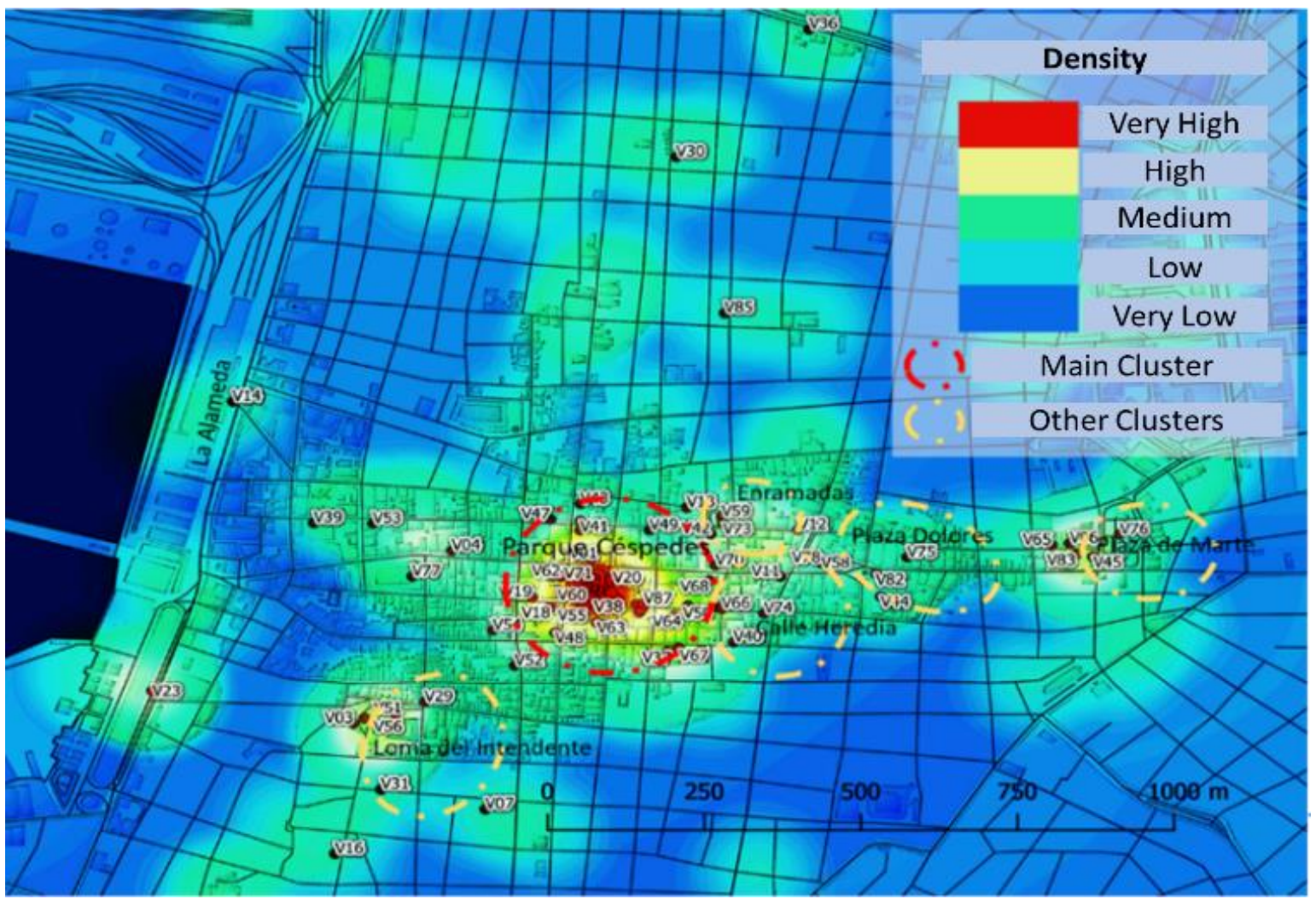

**Figure 4**. Heatmap of attraction concentration levels

### *4.3. Discussion*

The implementation of the proposed technique in Santiago de Cuba's tourism sector appears to contribute to the identification and prioritization of tourist attractions, which could be valuable for a developing tourist destination such as Santiago de Cuba. The FTV index not only revealed the experiential tourism value of the attractions in this destination but also made it possible to rank and filter, from the total set of attractions, those with greater potential for designing sustainable tourism experiences. Then, tourism managers can leverage attractions with high performance in the analyzed factors, and work to enhance sites of lower tourism value, potentially transforming them into valuable additions to the main sustainable tourism experience.

The FTV index use in each factor revealed which dimensions achieve a higher score and, consequently, contribute more positively to the selection of attractions. The high positive evaluations of the variables "Cultural Commercialization", "Social Vulnerability", and "Environmental Impact" indicate that there are no visible damages or adverse effects on the tourism destination's attractions and the residents view tourism positively. Therefore, it is recommended to take care of these factors to improve the visitor's experience.

From the results, the "Visitor Type" factor also stands out, which suggests that the types of visitors arriving at the Santiago de Cuba destination have a greater influence on the perception of a positive experience, highlighting the need for future experiential marketing actions in the destination focused on specific segments. Similarly, it is inferred that the "Attraction Info" and "Attraction Condition" factors are key determinants, therefore, this destination must monitor multi-channel promotion and prioritize regular maintenance and restoration of attractions. Additionally, the "Visit Date" factor unexpectedly reached a high FTV, which could be explained by the seasonal variations of the Santiago de Cuba destination, which in turn suggest the need for strategies to be considered to manage crowds during "peak" periods where traditional events or festivals coincide in the same period (for example, the Caribbean Festival and the Carnivals, both traditionally held in July).

Moreover, by combining the calculation of the proposed FTV with a Geographic Information System (GIS), enables a spatial data analysis for improving tourism planning. We can to optimize the use of tourism potential of different areas and designate special development zones for prioritized investment. Stakeholders can use these data to design or refine sustainable tourism experiences, making them more appealing and unique.

## 5. Conclusions and future work

The method proposed for building a Fuzzy Experiential Tourism Value index consists on a combination of a Triangular Fuzzy Rescaling (TFR) and a weighted average aggregation of Triangular Fuzzy Numbers. This approach has two main advantages in comparison with the usual crisp Tourism Index Values. On the one hand, it enables the consideration of subjective and imprecise factors, because we can either associate fuzzy numbers to the numerical scale used by the respondent, or we can give directly the opportunity to use a fuzzy linguistic variable constructed using fuzzy numbers. On the other hand, it permits the use of different scales of measurement on the factors, which are later transformed to a common scale with the TFR operation. Therefore, FTV can perform well in numerous contexts where diverse information exists, both objective and subjective, as is the case with sustainable tourism experiences. This makes it a powerful tool for measuring the level of tourist experiences in different tourist destinations and for establishing comparisons between them, regardless of differences in measurement scales or the degree of objectivity in each evaluation. This is because the FTV always aggregates the standardized and rescaled results on a dimensionless scale.

By involving a data normalization process, it is not necessary for these factors to be formulated in the same way when deciding to generalize the study to another destination. In other words, the proposed technique not only offers greater flexibility but also ensures a higher capacity for generalization. This could have important practical applications, from a cost-effectiveness perspective. Once its effective implementation is understood, it contributes to saving time and resources that are usually dedicated to the construction of different methodologies adjusted to each tourist destination. Instead, it is possible to use the same technique in heterogeneous destinations. The paper has shown the application of the proposed technique in the destination of Santiago de Cuba allowed for validating its utility in the practice of tourism management.

This study is a first step towards the definition of new tourism value indices, and needs more research as it has some limitations. The presented work was mainly focused on Santiago de Cuba and therefore, further research is also suggested to compare different destinations to validate the practical effectiveness of the proposed technique and to propose methodologies that can complement our approach for a more complete design. Moreover, for practical purposes, the proposed FTV simply uses a weighted arithmetic mean as an aggregation operator, which is compensative in nature. This means that low values in some factors can be offset by high values in others. However, in certain scenarios, decision-makers may prefer to employ non-compensative operators that do not allow trade-offs between criteria. Future research could explore the application of non-compensative aggregation operators, such as the Graded Conjunctive-Disjunctive (GCD) or the Ordered Weighted Averaging (OWA) operator.

**Acknowledgements**

This research received support from the *Agència de Gestió d'Ajuts Universitaris i de Recerca* (AGAUR) and was co-funded by the European Union, under Grant 2023 FI-1 00622, and takes part of the Grant PID2021-122575NB-I00 funded by MCIN/AEI/ 10.13039/501100011033/ by "ERDF A way of making Europe". This work was also partially funded by Universitat Rovira i Virgili (URV) [projects number 2023PFR-URV-114 and 2022PFR-URV-41], and by Generalitat de Catalunya (AGAUR) to ITAKA consolidated research group [2021-SGR-00114]. The support of the Spanish Research Network on Recommender Systems (ELIGE-IA, RED2022-134302-T) is also acknowledged.